\documentstyle[twoside,amsart12]{amsart}

\newtheorem{thm}{Theorem}[section]
\newtheorem{cor}[thm]{Corollary}
\newtheorem{lma}[thm]{Lemma}

\theoremstyle{definition}
\newtheorem{defn}[thm]{Definition}
\newtheorem{eg}[thm]{Example}
\newcommand\Hopf{{\mathop{\mathrm {\hbox{-}Hopf}}}} 
\newcommand\Lie{{\mathop{\mathrm {\hbox{-}Lie}}}} 
\newcommand\Der{{\mathop{\mathrm {Der}}}} 
\newcommand\End{{\mathop{\mathrm {End}}}} 
\newcommand\Alg{{\mathop{\mathrm {\hbox{-}Alg}}}} 
\newcommand\sgn{{\mathop{\mathrm {sgn}}}} 
\newcommand\tensor{\otimes} 
 
\newcommand\M{{\cal M}} 
\newcommand\N{{\Bbb N}} 
\newcommand\Z{{\Bbb Z}} 
 
\newcommand\z{{\mbox{-}}} 
\newcommand\x{{\overline x}} 
 
\def\=>{\longrightarrow}

\begin{document} 

\title[On Lie Algebras in Braided Categories]{On Lie 
Algebras in Braided Categories} 
\author{Bodo Pareigis} 
\address{Mathematisches Institut der Universit\"at 
M\"unchen\\ Germany} 
\email{pareigis@@rz.mathematik.uni-muenchen.de} 
\maketitle 

 
 \begin{abstract} 
 The category of group-graded modules over an abelian group 
$G$ is a monoidal category. For any bicharacter of $G$ this 
category becomes a braided monoidal category. We define the 
notion of a Lie algebra in this category generalizing the 
concepts of Lie super and Lie color algebras. Our Lie 
algebras have $n$-ary multiplications between various 
graded components. They possess universal enveloping 
algebras that are Hopf algebras in the given category. 
Their biproducts with the group ring are noncommutative 
noncocommutative Hopf algebras some of them known in the 
literature. Conversely the primitive elements of a Hopf 
algebra in the category form a Lie algebra in the above 
sense. 
 
 \end{abstract} 

\section{Introduction} 

 With the appearance of quantum groups the study of Hopf 
algebras has taken an important turn in recent years. Many 
families of quantum groups are known. Here we will add some 
new families by means of a construction which might 
eventually also help to develop a structure theory of 
quantum groups. In classical structure theory a formal 
group is decomposed into a smash product of an 
infinitesimal part with a separable part. The separable 
part is often a group algebra, whereas the infinitesimal 
part is often a universal enveloping algebra on which the 
group algebra operates. The infinitesimal part then is 
generated by the primitive elements of the formal group. 
 
 However, for quantum groups or noncommutative 
noncocommutative Hopf algebras this kind of decomposition 
seems to be much more complicated. We will use the 
following approach. 
 
 Let $\chi: G \times G \longrightarrow k^*$ be a 
bicharacter of the abelian group $G$. Then the monoidal 
category of $G$-graded vector spaces or $kG$-comodules is a 
braided monoidal category. By forming a biproduct with the 
group algebra $kG$, each Hopf algebra in this category 
generates an ordinary Hopf algebra. This process was 
studied for Hopf algebras with a projection by Radford in 
\cite{R} (or as the process of bosonization by Majid in \cite{M94a}). 
 
 We will give some general techniques for constructing Hopf 
algebras in the category of $G$-graded vector spaces. Our 
special interest lies in Hopf algebras that are generated 
by primitive elements in the proper sense, that is 
$\Delta(x) = x \tensor 1 + 1 \tensor x$. In the biproduct 
these elements become skew primitive elements. 
 
 There is a generalized notion of a Lie algebra in 
symmetric monoidal categories under the name of Lie color 
algebra (or Lie super algebra). An appropriate definition 
of a Lie algebra in a {\it braided} monoidal category 
should have specific properties. In the classical situation 
the primitive elements of a Hopf algebra form a Lie 
algebra. Also the set of derivations of an algebra is a Lie 
algebra. So a more general definition of a Lie algebra 
should pass a test with respect to the set primitive 
elements of a Hopf algebra and with respect to the set of 
derivations of an algebra. We propose a definition of 
generalized Lie algebras in the  braided monoidal category 
of $G$-graded vector spaces. These generalized Lie algebras 
have $n$-ary bracket multiplications that are only 
partially defined, but certain symmetry and Jacobi 
identities still hold. 
 
 Lie color and Lie super algebras as well as ordinary Lie 
algebras are special cases of our generalized Lie algebras. 
We will show that the set of primitive elements of a Hopf 
algebra is such a generalized Lie algebra. Every 
associative algebra is also a generalized Lie algebra by 
the same definition of the bracket multiplication. And the 
set of derivations of an algebra will turn out to be a Lie 
algebra. 
 
 Starting out with a generalized Lie algebra we will 
construct a universal enveloping algebra that turns out to 
be a Hopf algebra in the category of $G$-graded vector 
spaces. Thus we have obtained a method for constructing new 
Hopf algebras in this category. By forming biproducts we 
obtain many old and new ordinary noncommutative 
noncocommutative Hopf algebras or quantum groups. 
 
 Ordinary Hopf algebras over a field of characteristic zero 
that are generated by their primitive elements are either 
trivial or infinite dimensional. In fact they are universal 
enveloping algebras of Lie algebras. 
 
 Hopf algebras in braided monoidal categories in 
characteristic zero that are generated by their primitive 
elements, however, do not share this property. In certain 
cases they are similar to universal restricted enveloping 
algebras of restricted (p-)Lie algebras. 
 
 A simple example of such a Hopf algebra is ${\bf 
C}[x]/(x^n)$ with $\Delta(x) = x \tensor 1 + 1 \tensor x$ 
considered as an object in the braided monoidal category of 
$C_n$-graded vector spaces, where $C_n$ is the cyclic group 
of order $n$ and the braiding is given by a primitive 
$n$-th root of unity. It is generated by the primitive 
element $x$ as an algebra in the given category, but it is 
finite dimensional. Such a Hopf algebra cannot exist as an 
ordinary Hopf algebra (in the category of vector spaces). 
 
 We will restrict our considerations to group graded vector 
spaces over a fixed abelian group $G$. A consequence of one 
of our main results is that for primitive elements $x_1, 
\ldots, x_n$ of degree 1 in a $C_n$-graded Hopf algebra $H$ 
the expression 
 $$[x_1, \ldots, x_n] := \sum_{\sigma \in S_n} x_{\sigma(1)} 
\ldots x_{\sigma(n)}$$ 
 is again a primitive element (of degree zero) in $H$ 
(where $S_n$ is the symmetric group). 
 
 By now the reader should be interested to see the 
definition of a generalized $G$-graded Lie algebra. Let $G$ 
be an abelian group with a bicharacter $\chi: G \times G 
\longrightarrow k^*$. For every primitive $n$-th root of 
unity $\zeta$ we consider certain $n$-tuples $(g_1, \ldots, 
g_n)$ in $G$ associated with $\zeta$. We will call them 
$\zeta$-families. Furthermore we construct factors 
$\rho(\sigma,(g_1, \ldots, g_n)) \in k^*$ for each 
permutation $\sigma \in S_n$ and each $\zeta$-family $(g_1, 
\ldots, g_n)$. These factors generalize the sign of a 
permutation. Details will be given in Definition \ref{rho}. A 
generalized $G$-graded Lie algebra is a $G$-graded vector 
space $P  = \bigoplus_{g \in G} P_g$ that has multilinear 
bracket operations for all $\zeta$ and all $\zeta$-families 
$(g_1, \ldots, g_n)$ 
 $$[ \ldots ]: P_{g_1} \times \ldots \times P_{g_n} 
\longrightarrow P_{g_1 + \ldots + g_n}$$
 satisfying the following generalizations of the symmetry 
and Jacobi identities 
 \begin{itemize}
 \item $[x_1, \ldots, x_n] = \rho(\sigma,(g_1, 
\ldots, g_n))\ [x_{\sigma(1)}, \ldots, x_{\sigma(n)}],$ 
 \item $\sum_{i=1}^{n+1} \zeta^{-i+1} \big( 
\prod_{j=1}^{i-1} \chi(g_j,g_i) \big)\ 
[x_i,[x_1,\ldots,\hat x_i, \ldots, x_{n+1}]] = 0,$ 
 \item $[x,[y_1, y_2, \ldots, y_n]] = 
\sum_{i=1}^n \big( \prod_{j=1}^{i-1} \chi(g_j,h) \big)\ 
[y_1, \ldots, [x,y_i], \ldots, y_n],$
 \end{itemize}
 \parindent0pt 
 
 whenever the bracket products are defined. In particular 
we require $x_i,y_i \in P_{g_i}$ and $x \in P_h$. 
\parindent15pt 
 
 The set of primitive elements of a $G$-graded Hopf algebra 
is an example for a generalized Lie algebra, and so is any 
$G$-graded algebra if we define the bracket operation by 
 $$[x_1, \ldots, x_n] := \sum_{\sigma \in S_n} \rho(\sigma, 
(g_1, \ldots, g_n)) x_{\sigma(1)} \cdot \ldots \cdot 
x_{\sigma(n)}.$$ 
 
 A Lie super algebra $(P_0,P_1)$, where $P_0$ is an 
ordinary Lie algebra with operations $[.,.]: P_1 \tensor 
P_1 \longrightarrow P_0$ and $[.,.]: P_0 \tensor P_1 
\longrightarrow P_1$, is an example for this definition. 
Lie color algebras are also special cases of our concept. 

 The theory of braided Lie algebras will be generalized to 
the category of Yetter-Drinfeld modules in \cite{P}.
 
 I would like to acknowledge helpful conversations with 
Edward L. Green, Frank Halanke, Susan Montgomery, Martin 
Neuchl, Helmut Rohrl, Peter Schauenburg, and Yorck 
Sommerh\"auser. 
 
 \section{Bicharacters and the Bracket Multiplication} 
 
 Throughout let $k$ be an integral domain, let $G$ be an 
abelian group written additively and let $\chi: G \times G 
\longrightarrow k^*$ be a bicharacter, that is a group 
homomorphism $\chi: G \tensor_{\bf Z} G \longrightarrow 
k^*$. Then bicharacters on $\Z$ resp. $\Z/n\Z$ have the 
form $\chi(i,j) = \zeta^{ij}$. 
 
 \begin{defn} \label{rho} 
 Let $\zeta \in k^*$. An $n$-tuple $(g_1, \ldots, g_n)$ of 
(not necessarily distinct) elements $g_i \in G$ is a 
$\zeta${\it -family} of {\it length} $n$, if 
 $$\chi(g_i,g_j)\chi(g_j,g_i) = \zeta^2$$ 
 for all $i \not= j$.
 \end{defn} 
 
 Observe that a $\zeta$-family $(g_1, \ldots, g_n)$ is also 
a $(-\zeta)$-family. 
 
 If $(g_1, \ldots, g_n)$ is a $\zeta$-family and $\sigma 
\in S_n$ then $\sigma(g_1, \ldots, g_n) := (g_{\sigma(1)}, 
\ldots, g_{\sigma(n)})$ is also a $\zeta$-family. Thus the 
set $G_n^\zeta$ of all $\zeta$-families of length $n$ is an 
$S_n$-set. Observe that there are $\zeta$-families of 
varying lengths $n$. 
 
 Let $X \subset G_n^\zeta$ be an $S_n$-closed subset. We 
define a map $\rho: S_n \times X \longrightarrow k^*$ by 
 $$\rho(\sigma,(g_1, \ldots, g_n)) := \prod_{(i,j) \in R} 
(\zeta^{-1} \chi(g_{\sigma(j)}, g_{\sigma(i)}))$$ 
 where $R$ is the set of all pairs $(i,j)$ with $1 \leq i < 
j \leq n$ and $\sigma(i) > \sigma(j)$, i.e. the set of all 
pairs in $\sigma$ in reverse position. 
 
 \begin{lma}\label{rho_is_homo} 
 The map $\rho$ satisfies the following relation 
 $$\rho(\sigma\tau,(g_1, \ldots, g_n)) 
= \rho(\tau,(g_{\sigma(1)}, \ldots, g_{\sigma(n)})) 
\rho(\sigma,(g_1, \ldots, g_n)).$$
 \end{lma} 

 \begin{pf} 
 We have to show 
 $$ \begin{array}{rl} 
 &\prod_{i<j, \sigma\tau(i)>\sigma\tau(j)} 
  (\zeta^{-1}\chi(g_{\sigma\tau(j)}, g_{\sigma\tau(i)})) \\ 
 &= \displaystyle\Big( \prod_{k<l,\tau(k)>\tau(l)} 
  (\zeta^{-1}\chi(g_{\sigma\tau(l)}, g_{\sigma\tau(k)})) 
\Big) 
 \Big( \prod_{r<s,\sigma(r)>\sigma(s)} 
  (\zeta^{-1}\chi(g_{\sigma(s)}, g_{\sigma(r)}))\Big).\\ 
 \end{array}
 $$ 
 We investigate the factors on the left hand side. If $i<j$ 
with $\sigma\tau(i) < \sigma\tau(j)$ then there is no 
corresponding factor on the l.h.s. On the r.h.s. there are 
now two possibilities. If $s:= \tau(i) > \tau(j) =: r$ (and 
thus $\sigma(r) > \sigma(s)$) then the first product 
contributes a factor $\zeta^{-1} \chi(g_{\sigma\tau(j)}, 
g_{\sigma\tau(i)})$ for the pair $i<j, \tau(i) > \tau(j)$ 
and the second product contributes a factor $\zeta^{-1} 
\chi(g_{\sigma\tau(i)},g_{\sigma\tau(j)})$ for the pair 
$r<s, \sigma(r)>\sigma(s)$. These two factors cancel. If 
$s:= \tau(i) < \tau(j) =: r$ then neither pair $i<j, 
\tau(i) < \tau(j)$ nor $r<s, \sigma(r) < \sigma(s)$ 
contributes a factor. 
 
 However, if $i<j$ with $\sigma\tau(i) > \sigma\tau(j)$ 
then there is a corresponding factor $\zeta^{-1} 
\chi(g_{\sigma\tau(j)}, g_{\sigma\tau(i)})$ on the l.h.s. 
On the r.h.s. there are again two possibilities. If 
$\tau(i) > \tau(j)$ then the first product contributes a 
factor $\zeta^{-1} \chi(g_{\sigma\tau(j)}, 
g_{\sigma\tau(i)})$ for the pair $i<j, \tau(i) > \tau(j)$. 
The second product does not contribute a factor for the 
pair $s:= \tau(i) > \tau(j) =: r, \sigma(r) < \sigma(s)$. 
If $s:= \tau(i) < \tau(j) =: r$ then the first product does 
not contribute a factor for the pair $i<j, \tau(i) < 
\tau(j)$, but the second product contributes a factor 
$\zeta^{-1} \chi(g_{\sigma\tau(i)},g_{\sigma\tau(j)})$ for 
the pair $r<s, \sigma(r)>\sigma(s)$. This argument takes 
care of all factors on both sides of the formula. 
 \end{pf} 
 
 If $g,g' \in G$ have order $m$ resp. $m'$ and if $n = 
\gcd(m,m')$ then $\chi(g,g')^n = 1$ since $\chi(g,g')^m = 
\chi(mg,g') = \chi(0,g') = 1$ and $\chi(g,g')^{m'} = 
\chi(g,m'g') = 1$. Hence $\chi(g,g')$ is a primitive $n$-th 
root of unity for all $g,g' \in G$ of finite order, with 
$n$ chosen suitably. 
 
 Apart from this restriction on the choice of $\chi(g,g')$ 
any combination of values can occur. In particular let 
$\zeta$ be an $n$-th root of unity. Then there are examples 
of groups $G$ and elements $g_1, \dots, g_m \in G$ with 
$\chi(g_i,g_j) = \zeta$ for all $i \not= j$. Take for 
example $G = C_n \oplus \ldots \oplus C_n = \Z g_1 \oplus 
\ldots \oplus \Z g_m$. Then a bicharacter $\chi$ is defined 
as the group homomorphism $\chi: G \tensor G = 
\bigoplus_{i,j = 1}^m \Z g_i \tensor \Z g_j = 
\bigoplus_{i,j = 1}^m \Z(g_i \tensor g_j) \longrightarrow 
k^*$, where the images $\chi(g_i \tensor g_j) = 
\chi(g_i,g_j)$ are chosen arbitrarily among the $n$-th 
roots of unity (including 1). 
 
 If $(g_1, \ldots, g_n)$ is a $\zeta$-family with at least 
one element $g_i$ of finite order, then 
$\chi(g_i,g_j)\chi(g_j,g_i) = \zeta^2$ implies that $\zeta$ 
is a root of unity. 
 
 \begin{eg} 
 For $g \in G$ we define $|g| := \chi(g,g)$. If $g$ has 
order $m$ then $|g|$ is a primitive $n$-th root of unity 
where $n$ divides the order $m$ of $g$. Furthermore 
$(g,\ldots,g)$ is a $|g|$-family in $G$. Observe also that 
any pair $(0,g)$ or $(g,0)$ is a 1-family and a 
 $(-1)$-family. 
 
 The map $\rho: S_n \times X \longrightarrow k^*$ reduces 
to a well known map in the following situation. If 
$(g_1,\ldots,g_n) = (g,\ldots,g)$, $|g| = 1$, and $\zeta = 
-1$ then $\rho(\sigma, (g,\ldots,g)) = \sgn(\sigma)$ is the 
sign of the permutation. 
 \end{eg} 
 
 Now we turn to the category of interest for us. By 
\cite{FM} Remark 3.4 the category $\M^{kG}$ of $G$-graded 
$k$-modules (or the category of $kG$-comodules) is a 
braided monoidal category with the tensor product 
 $$(X \tensor Y)_g = \bigoplus_{h \in G}X_h \tensor 
Y_{g-h}$$ 
 and the braiding 
 $$\tau: X \tensor Y \ni x \tensor y \mapsto 
\chi(\deg(x),\deg(y)) y \tensor x \in Y \tensor X$$ 
 where $x$ and $y$ are homogeneous elements of degree 
$\deg(x), \deg(y) \in G$. 
 
 Algebras, coalgebras, bialgebras, and Hopf algebras in the 
braided monoidal category $\M^{kG}$ will be called 
$(G,\chi)$-algebras, $(G,\chi)$-coalgebras, 
$(G,\chi)$-bialgebras, resp. $(G,\chi)$-Hopf algebras. 
 
 Let $A$, $B$ be $(G,\chi)$-algebras in the category 
$\M^{kG}$. Then $A \tensor B$ is also a $(G,\chi)$-algebra 
in $\M^{kG}$ with the multiplication $A \tensor B \tensor A 
\tensor B \buildrel 1 \tensor \tau \tensor 1 \over 
\longrightarrow A \tensor A \tensor B \tensor B \buildrel 
m_A \tensor m_B \over \longrightarrow A \tensor B$ (see for 
example \cite{M94b}). 
 
 \begin{defn} \label{assLie} 
 Let $A$ be a $(G,\chi)$-algebra (associative with unit) in 
$\M^{kG}$. We define 
 $$[x_1, \ldots, x_n] := \sum_{\sigma \in S_n} \rho(\sigma, 
(g_1, \ldots, g_n)) x_{\sigma(1)} \cdot \ldots \cdot 
x_{\sigma(n)}$$ 
 for all $\zeta \in k^*$, all 
$\zeta$-families $(g_1, \ldots, g_n)$, and all $x_i \in 
A_{g_i}$. 
 \end{defn} 
 
 A special example is $\zeta = -1$, $n = 2$ and 
$\chi(g_1,g_2) = \chi(g_2,g_1) = 1$. Then 
 $$[x_1,x_2] = x_1x_2 + \rho((2,1),(g_1,g_2))x_2x_1 = 
x_1x_2 - x_2x_1.$$ 
 If $\zeta = -1$, $n = 2$ and $\chi(g_1,g_2) = 
\chi(g_2,g_1) = -1$, then 
 $$[x_1,x_2] = x_1x_2 + \rho((2,1),(g_1,g_2))x_2x_1 = 
x_1x_2 + x_2x_1.$$ 
 Hence we obtain Lie algebras resp. Lie super algebras as a 
special case of (\ref{assLie}). Similarly if $\zeta = -1$, $n = 2$ 
and $(g_1,g_2)$ is a $\zeta$-family then $\chi(g_1,g_2) = 
\chi(g_2,g_1)^{-1}$ and 
 $$[x_1,x_2] = x_1x_2 - \chi(g_1,g_2) x_2x_1,$$ 
 rendering Lie color algebras as another instance of (\ref{assLie}). 
For Lie color algebras the assumption $\chi(g_2,g_1) = 
\chi(g_1,g_2)^{-1}$ for all $g_1,g_2 \in G$ leads to a 
totally defined multiplication $[.,.]: A \times A 
\longrightarrow A$ since $[x_1,x_2]$ is defined for any 
pair $x_1,x_2 \in A$ of homogeneous elements. A special new 
operation is the following bracket product. If $g \in G$ 
satisfies $\chi(g,g) = \zeta \not= 1$ then we have the 
product $[ \ldots ]: \tensor^n A \longrightarrow A$ 
 $$[x_1, \ldots, x_n] = \sum_{\sigma \in S_n} x_{\sigma(1)} 
\ldots x_{\sigma(n)}.$$ 
 
 \begin{thm}(Symmetry)\label{symm} 
 Let $\zeta \in k^*$ and $(g_1, \ldots, g_n)$ be a 
$\zeta$-family. Let $x_i \in A_{g_i}$ and $\sigma \in 
S_n$. Then 
 $$[x_1, \ldots, x_n] = \rho(\sigma,(g_1, \ldots, g_n)) 
[x_{\sigma(1)}, \ldots, x_{\sigma(n)}].$$
 \end{thm} 
 
 \begin{pf}
 We have 
 $$\begin{array}{rl} 
 \rho(\sigma,(g_1, \ldots, g_n)) &\hskip-1ex 
[x_{\sigma(1)}, \ldots, x_{\sigma(n)}] =\\ 
 &\hskip-2ex = \sum_{\tau \in S_n} 
 \rho(\sigma,(g_1, \ldots, g_n)) \rho(\tau, (g_{\sigma(1)}, 
\ldots, g_{\sigma(n)})) 
 x_{\sigma\tau(1)} \cdot \ldots \cdot x_{\sigma\tau(n)} \\ 
 &= \sum_{\tau \in S_n} 
 \rho(\sigma\tau,(g_1, \ldots, g_n)) 
 x_{\sigma\tau(1)} \cdot \ldots \cdot x_{\sigma\tau(n)}. 
 \qed 
 \end{array}$$ 
 \renewcommand{\qed}{} \end{pf} 
 
 In the case of Lie algebras this amounts to 
 $$[x_1,x_2] = -[x_2,x_1];$$ 
 in the case of Lie super algebras this is 
 $$[x_1,x_2] = [x_2,x_1]$$ 
 (for $\zeta = -1$, $n = 2$ and $\chi(g_1,g_2) = 
\chi(g_2,g_1) = -1$); and in the case of Lie color algebras 
(see \cite {FM} 3.11) this is 
 $$[x_1,x_2] = -\chi(g_1,g_2)[x_2,x_1].$$ 
 
 In the following theorem let $(i \ldots 1) = {{1, 2, 
\ldots, i, i+1, \ldots, n+1} \choose {i, 1, \ldots, i-1, 
i+1, \ldots, n+1}}$ denote a cycle in $S_{n+1}$. 

 \begin{thm}{(Jacobi identities)} \label{jacobi}
 \begin{enumerate} 
 \item {1.} 
 Let $(g_1, \ldots, g_{n+1})$ be a $\zeta$-family with 
$\zeta$ a primitive $n$-th root of unity. Then 
 $$\sum_{i=1}^{n+1} \rho((i \ldots 1), (g_1, \ldots, 
g_{n+1})) [x_i,[x_1,\ldots,\hat x_i, \ldots, x_{n+1}]] = 
0$$ 
 for all $x_i \in A_{g_i}$. 
 \item{2.} 
 Let $(g_1, \ldots, g_n)$ be a $\zeta$-family with $\zeta$ 
a primitive $n$-th root of unity and let $h \in G$ such 
that all $(h,g_i)$ are $(-1)$-families. Then 
 $$[x,[y_1, y_2, \ldots, y_n]] = \sum_{i=1}^n \big( 
\prod_{j=1}^{i-1} \chi(g_j,h) \big) [y_1, \ldots, [x,y_i], 
\ldots, y_n] $$ 
 \item{} for all $x \in A_h$ and $y_i \in A_{g_i}$. 
 \end{enumerate}
 \end{thm}
 
 \begin{pf} 
 (1) Let $\tau \in S_n$. Construct $\bar\tau \in S_{n+1}$ 
by $\bar\tau(1) = 1$ and $\bar\tau(j) = \tau(j-1)+1$. This 
defines a bijection between $S_n$ and the set of all 
$\sigma \in S_{n+1}$ with $\sigma(1) = 1$. Using $h_i := 
g_{i+1}$ we get 
 $$\begin{array}{rl} 
 \rho(\tau, (g_2, \ldots, g_{n+1})) &= \rho(\tau, (h_1, 
\ldots, h_n)) \\ 
 &= \displaystyle \prod_{1 \leq i < j \leq n, \tau(i) > 
\tau(j)} (\zeta^{-1} \chi(h_{\tau(j)},h_{\tau(i)})) \\ 
 &= \displaystyle \prod_{2 \leq i < j \leq n+1, \tau(i-1) > 
\tau(j-1)} (\zeta^{-1} \chi(g_{\tau(j-1)+1},g_{\tau(i- 
1)+1})) \\ 
 &= \displaystyle \prod_{1 \leq i < j \leq n+1, \bar\tau(i) 
> \bar\tau(j)} (\zeta^{-1} 
\chi(g_{\bar\tau(j)},g_{\bar\tau(i)})) \\ 
 &= \rho(\bar\tau, (g_1, \ldots, g_{n+1})) \\
 \end{array} 
 $$ 
 So we get for all $i \in \{1, \ldots, n+1\}$ 
 $$\rho(\tau, (g_1, \ldots, g_{i-1}, \hat g_i, g_{i+1}, 
\ldots, g_{n+1})) = \rho(\bar\tau, (g_i, g_1, \dots, g_{i- 
1}, \hat g_i, g_{i+1}, \ldots, g_{n+1})).$$ 
 Now let $\sigma = \sigma_i \in S_{n+1}$ with $i := 
\sigma(1)$. Let $\bar\tau := (1 \ldots i)\sigma_i$. Then 
$\bar\tau(1) = 1$, so $\bar\tau$ comes from some $\tau \in 
S_n$. Furthermore $\sigma_i = (i \ldots 1)\bar\tau$. By 
Lemma \ref{rho_is_homo} we get 
 $$\begin{array}{rl} 
  \rho(\sigma_i,(g_1, \ldots, g_{n+1})) 
 &= \rho(\bar\tau,(g_i, g_1, \dots, \hat g_i, \ldots, 
g_{n+1})) 
  \rho((i \ldots 1),(g_1, \ldots, g_{n+1}))  \\ 
 &= \rho(\tau, (g_1, \ldots, \hat g_i, \ldots, g_{n+1})) 
  \rho((i \ldots 1),(g_1, \ldots, g_{n+1})). \\
 \end{array}
 $$ 
 
 Given $\tau \in S_n$ we define $\tilde\tau$ by 
$\tilde\tau(n+1) := n+1$ and $\tilde\tau(j) := \tau(j)$ 
else. This defines a bijection between $S_n$ and the set of 
all $\sigma \in S_{n+1}$ with $\sigma(n+1) = n+1$. Then we 
get 
 $$\rho(\tau,(g_1, \ldots, g_n)) = \rho(\tilde\tau, (g_1, 
\ldots, g_{n+1}))$$ and 
 $$\rho(\tau, (g_1, \ldots, \hat g_i, \ldots, g_{n+1})) = 
 \rho(\tilde\tau, (g_1, \dots, \hat g_i, \ldots, g_{n+1}, 
g_i)).$$ 
 Now let $\sigma = \sigma^i \in S_{n+1}$ with $i := 
\sigma(n+1)$. Let $\tilde\tau := (n+1 \ldots i)\sigma^i$. 
Then $\tilde\tau(n+1) = n+1$, so $\tilde\tau$ comes from 
some $\tau \in S_n$. Furthermore $\sigma^i = (i \ldots 
n+1)\tilde\tau$ and we get 
 $$ \rho(\sigma^i,(g_1, \ldots, g_{n+1})) 
 = \rho(\tau, (g_1, \ldots, \hat g_i, \ldots, g_{n+1})) 
\rho((i \ldots n+1),(g_1, \ldots, g_{n+1})).$$ 
 
 To prove the Jacobi identity we observe that since $\zeta$ 
is a primitive $n$-th root of unity 
 $$\begin{array}{r@{}l}
 \chi(g_i,g_1 + \ldots + \hat g_i + \ldots & + g_{n+1}) 
 \chi(g_1 + \ldots + \hat g_i + \ldots + g_{n+1},g_i) \\ 
 &= \prod_{j \not= i}\chi(g_i,g_j) \chi(g_j,g_i) = 
\zeta^{2n} = 1, \\
 \end{array}
 $$ 
 hence we have a $(-1)$-family $(g_i,g_1 + \ldots + \hat 
g_i + \ldots + g_{n+1})$ and 
 $$\begin{array}{rl} [x_i,&[x_1, \ldots, \hat x_i, \ldots, 
x_{n+1}]] = \\ 
 &= 
 x_i[x_1, \ldots, \hat x_i, \ldots, x_{n+1}] - \big( 
\prod_{j \not= i}\chi(g_i,g_j) \big) 
 [x_1, \ldots, \hat x_i, \ldots, x_{n+1}]x_i. \\
 
\end{array}
 $$ 
 Furthermore we use $\zeta^n = 1$ to get 
 $$\begin{array}{rl}
 \rho((i \ldots 1),(g_1, \ldots, g_{n+1})) 
 &\prod_{j \not= i} \chi(g_i,g_j) = \\ 
 &= \big( \prod_{j<i} (\zeta^{-1} \chi(g_j,g_i)) \big) 
 \big( \prod_{j \not= i} (\zeta^{-1} \chi(g_i,g_j)) \big) 
\\ 
 &= \prod_{j>i} (\zeta^{-1} \chi(g_i,g_j)) = \rho((i \ldots 
n+1),(g_1, \ldots, g_{n+1})). \\
 \end{array}
 $$ 
 With these notations we can now evaluate 
 $$ \begin{array}{l}
 \sum_{i=1}^{n+1} \rho((i \ldots 1),(g_1, \ldots, 
g_{n+1})) 
   [x_i,[x_1, \ldots, \hat x_i, \ldots, x_{n+1}]] = \\ 
 = \sum_{i=1}^{n+1} 
    \rho((i \ldots 1),(g_1, \ldots, g_{n+1})) 
    x_i[x_1, \ldots, \hat x_i, \ldots, x_{n+1}]  \\ 
 \ \ - \rho((i \ldots 1),(g_1, \ldots, g_{n+1})) 
    \big( \prod_{j \not= i}\chi(g_i,g_j) \big) 
    [x_1, \ldots, \hat x_i, \ldots, x_{n+1}]x_i \\ 
 = \sum_{i=1}^{n+1} \sum_{\tau \in S_n} 
    \rho((i \ldots 1),(g_1, \ldots, g_{n+1})) 
    \rho(\tau,(g_1, \ldots, \hat g_i, \ldots, g_{n+1})) 
    x_{\sigma_i(1)}  \ldots 
    x_{\sigma_i(n+1)} \\ 
 \ \ - \rho((i \ldots n+1),(g_1, \ldots, g_{n+1})) 
    \rho(\tau,(g_1, \ldots, \hat g_i, \ldots, g_{n+1})) 
    x_{\sigma^i(1)} \ldots x_{\sigma^i(n+1)} \\ 
 = \sum_{\sigma \in S_{n+1}} \rho(\sigma,(g_1, \ldots, 
g_{n+1})) 
    x_{\sigma(1)} \ldots x_{\sigma(n+1)} 
 \ - \rho(\sigma,(g_1, \ldots, g_{n+1})) 
    x_{\sigma(1)} \ldots x_{\sigma(n+1)} \\ 
 = 0. \\
 \end{array}
 $$ 
 
 (2) Since $\chi(h,g_i)\chi(g_i,h) = (-1)^2$ and 
$\chi(g_i,g_j) \chi(g_j,g_i) = \zeta^2$ we have $\chi(h,g_1 
+ \ldots + g_n) \chi(g_1 + \ldots + g_n,h) = 1$ and 
$\chi(g_i,h + g_j)\chi(h + g_j,g_i) = \zeta^2$ so that all 
terms are defined. 
 
 Let $\sigma \in S_n$ with $\sigma(j) = i$. Then we have 
 $$ \begin{array}{rl} 
  \rho(\sigma,&(g_1, \ldots, g_{i-1}, h+g_i, g_{i+1}, 
\ldots, g_n)) = \\ 
 &= \displaystyle \big( \prod_{k < l, \sigma(k) > 
\sigma(l)} 
 (\zeta^{-1} \chi ( g_{\sigma(l)}, g_{\sigma(k)} )) \big) 
   \big( \prod_{j < l, i > \sigma(l)} \chi( g_{\sigma(l)}, 
h ) \big) 
   \big( \prod_{k < j, \sigma(k) > i} \chi( h, 
g_{\sigma(k)} ) \big) \\ 
 &= \displaystyle \rho(\sigma,(g_1, \ldots, g_n)) 
   \big( \prod_{\sigma^{-1}(i) < l, i > \sigma(l)} \chi( 
g_{\sigma(l)}, h ) \big) 
   \big( \prod_{k < \sigma^{-1}(i), \sigma(k) > i} \chi( h, 
g_{\sigma(k)} ) \big). \\
 \end{array}
 $$ 
 We abbreviate $z_k := u_k := y_k$ for $k \not= i$, $z_i := 
xy_i$, and $u_i := y_ix$. Then we get 
 $$ \begin{array}{rl}
 &\sum_{i=1}^n  
 \big( \prod_{r=1}^{i-1} \chi(h, g_r) \big) 
    [y_1, \ldots, [x,y_i], \ldots, y_n] \\ 
 &= \displaystyle \sum_{i=1}^n \big( \prod_{r=1}^{i-1} 
\chi(h, g_r) \big) 
    [y_1, \ldots, xy_i - \chi( h, g_i)y_ix, \ldots, y_n] 
\\ 
 &= \displaystyle \sum_{i=1}^n \sum_{\sigma \in S_n} \Big( 
   \big( \prod_{r=1}^{i-1} \chi(h, g_r) \big) 
   \rho(\sigma, (g_1, \ldots, h+g_i, \ldots, g_n)) 
   z_{\sigma(1)} \ldots z_{\sigma(n)} \\ 
 &\ \ - \big( \prod_{r=1}^{i-1} \chi(h, g_r) \big) 
   \chi( h, g_i) \rho(\sigma, (g_1, \ldots, h+g_i, \ldots, 
g_n)) 
   u_{\sigma(1)} \ldots u_{\sigma(n)}\Big) \\ 
 \end{array}
 $$ 

 $$ \begin{array}{rl}
 &= \displaystyle \sum_{\sigma \in S_n} \sum_{j=1}^n 
   \big( \prod_{r=1}^{\sigma(j)-1} \chi(h, g_r) \big) 
   \rho(\sigma,(g_1, \ldots, g_n)) 
   \big( \prod_{j < l, \sigma(j) > \sigma(l)} \chi( 
g_{\sigma(l)}, h ) \big) \cdot \\ 
 & \qquad \cdot  \big( \prod_{l < j, \sigma(l) > \sigma(j)} 
\chi( h, g_{\sigma(l)} ) \big) 
   y_{\sigma(1)} \ldots x y_{\sigma(j)} \ldots 
y_{\sigma(n)} \\ 
 &\ \ - \displaystyle \sum_{\sigma \in S_n} \sum_{j=1}^n 
\big( \prod_{r=1}^{\sigma(j)-1} \chi(h,g_r) \big) 
   \chi( h, g_{\sigma(j)}) \rho(\sigma,(g_1, \ldots, g_n)) 
   \big( \prod_{j < l, \sigma(j) > \sigma(l)} \chi( 
g_{\sigma(l)}, h ) \big) \cdot \\ 
 &\qquad \cdot \big( \prod_{l < j, \sigma(l) > \sigma(j)} 
\chi( h, g_{\sigma(l)} ) \big) 
   y_{\sigma(1)} \ldots y_{\sigma(j)} x \ldots 
y_{\sigma(n)} \\ 
 &=  \sum_{\sigma \in S_n} \rho(\sigma,(g_1, \ldots, g_n)) 
   xy_{\sigma(1)} \ldots y_{\sigma(n)} \\ 
 &\ \ - \displaystyle \sum_{\sigma \in S_n} \big( 
\prod_{r=1}^n \chi(h, g_r) \big) 
   \rho(\sigma,(g_1, \ldots, g_n)) 
   y_{\sigma(1)} \ldots y_{\sigma(n)}x \\ 
 &= [x,[y_1, \ldots, y_n]]. \\
 \end{array}
 $$ 
 if we can show that the coefficients reduce appropriately. 
For the first term $(j = 1)$ of the first sum and the last 
term $(j = n)$ of the second sum this is easy to see. The 
other terms cancel for all $\sigma \in S_n$ and all $j$. In 
fact, let $q := \sigma(j)$ and $p = \sigma(j-1)$. We have 
to show 
 $$\begin{array}{rl}
 \big( \prod_{r=1}^{\sigma(j)-1} 
 &\displaystyle \chi(h, g_r) \big) 
   \big( \prod_{j < l, \sigma(j) > \sigma(l)} \chi( 
g_{\sigma(l)}, h ) \big) 
   \big( \prod_{l < j, \sigma(l) > \sigma(j)} \chi( h, 
g_{\sigma(l)}) \big) \\ 
 &= \displaystyle \big( \prod_{r=1}^{\sigma(j-1)-1} \chi(h, 
g_r) \big) 
   \chi( h, g_{\sigma(j-1)}) 
   \big( \prod_{j-1 < l, \sigma(j-1) > \sigma(l)} \chi( 
g_{\sigma(l)}, h ) \big) \cdot \\ 
 &\ \ \cdot 
   \big( \prod_{l < j-1, \sigma(l) > \sigma(j-1)} \chi( h, 
g_{\sigma(l)} ) \big) \\
 \end{array}
  $$ 
 We change parameters by $\tau = \sigma^{-1}$ with 
$\tau(p)+1 = \tau(q)$ and have to show 
 $$\begin{array}{rl} 
 \big( \prod_{r=1}^{q-1} 
 &\chi(h,g_r) \big) 
   \big( \prod_{q > l, \tau(q) < \tau(l)} \chi( g_l, h ) 
\big) 
   \big( \prod_{l > q, \tau(l) < \tau(q)} \chi( h, g_l) 
\big) \\ 
 &= \big( \prod_{r=1}^p \chi(h, g_r) \big) 
   \big( \prod_{p > l, \tau(p) < \tau(l)} \chi( g_l, h ) 
\big) 
   \big( \prod_{l > p, \tau(l) < \tau(p)} \chi( h, g_l ) 
\big). \\
 \end{array}
  $$ 
 This can be easily checked if one considers the cases $p < 
q$ and $p > q$ separately. 
 \end{pf}
 
 
 \section{Primitive Elements} 
 
 We come to the main technical theorem of this paper which 
has applications to primitive elements in Hopf algebras. 
 
 \begin{thm} \label{main}
 Let $A$ be a $(G,\chi)$-algebra in $\M^{kG}$. Then the 
following hold in $A \tensor A$ 
 $$[x_1 \tensor 1 + 1 \tensor x_1, \ldots, x_n \tensor 1 + 
1 \tensor x_n] = [x_1, \ldots, x_n]  \tensor 1 + 1 \tensor 
[x_1, \ldots, x_n] $$ 
 for all primitive $n$-th roots of unity $\zeta$, all 
$\zeta$-families $(g_1, \ldots, g_n)$, and all $x_i \in 
A_{g_i}$. 
 \end{thm}
 
 \begin{pf} 
 We have to evaluate 
 $$\begin{array}{r@{}l} 
 [x_1 \tensor 1 &+ 1 \tensor x_1, \ldots, x_n \tensor 1 + 1 
\tensor x_n] \\ 
 &= \sum_{\sigma \in S_n} \rho(\sigma, (g_1, \ldots, g_n)) 
(x_{\sigma(1)} \tensor 1 + 1 \tensor x_{\sigma(1)}) \ldots 
(x_{\sigma(n)} \tensor 1 + 1 \tensor x_{\sigma(n)}). \\
 \end{array}
 $$ 
 Observe that 
 $$\begin{array}{rl} 
 &(x_i \tensor 1)(x_j \tensor 1) = (x_ix_j \tensor 1), \\ 
 &(x_i \tensor 1)(1 \tensor x_j) = (x_i \tensor x_j), \\ 
 &(1 \tensor x_i)(1 \tensor x_j) = (1 \tensor x_ix_j), \\ 
 &(1 \tensor x_i)(x_j \tensor 1) = \chi(g_i,g_j)(x_j 
\tensor x_i). \\
 \end{array}
 $$ 
 We collect all terms of the form $c_{\sigma,i} \cdot 
x_{\sigma(1)} \ldots x_{\sigma(i)} \tensor x_{\sigma(i+1)} 
\ldots x_{\sigma(n)}$ with $\sigma \in S_n$ and want to 
show that they are zero for all $1 < i < n$. 
 
 By Theorem \ref{symm} $[x_1 \tensor 1 + 1 \tensor x_1, \ldots, 
x_n \tensor 1 + 1 \tensor x_n]$ and $[x_{\sigma(1)} \tensor 
1 + 1\tensor x_{\sigma(1)}, \ldots, x_{\sigma(n)} \tensor 1 
+ 1 \tensor x_{\sigma(n)}]$ have (generically) proportional 
terms, in particular $c_{\sigma,i} \cdot x_{\sigma(1)} 
\ldots x_{\sigma(i)} \tensor x_{\sigma(i+1)} \ldots 
x_{\sigma(n)} = \rho(\sigma,(g_1, \ldots, g_n)) c_{1,i} 
\cdot x_{\sigma(1)} \ldots x_{\sigma(i)} \tensor 
x_{\sigma(i+1)} \ldots x_{\sigma(n)}$. So we may just 
consider the case $\sigma = 1$ or the term $c_{1,i} \cdot 
x_1 \ldots x_i \tensor x_{i+1} \ldots x_n$. 
 
 The term $c_{1,n} \cdot x_1 \ldots x_n \tensor 1$ occurs 
only in the product $(x_1 \tensor 1 + 1 \tensor x_1)(x_2 
\tensor 1 + 1 \tensor x_2) \ldots (x_n \tensor 1 + 1 
\tensor x_n)$ with the factor $c_{1,n} = 1$. The same holds 
for the term $c_{1,1} \cdot 1 \tensor x_1 \ldots x_n$. Now 
we consider exclusively terms $c_{1,i} \cdot x_1 \ldots x_i 
\tensor x_{i+1} \ldots x_n$ with $0 < i < n$. 
 
 We study which terms of the expansion of 
 $$\sum_{\sigma \in S_n} \rho(\sigma, (g_1, \ldots, g_n)) 
(x_{\sigma(1)} \tensor 1 + 1 \tensor x_{\sigma(1)}) \ldots 
(x_{\sigma(n)} \tensor 1 + 1 \tensor x_{\sigma(n)})$$ 
 contribute to $c_{1,i} \cdot x_1 \ldots x_i \tensor 
x_{i+1} \ldots x_n$. This will be those products of factors 
$x_1 \tensor 1, \ldots , x_i \tensor 1$ and $1 \tensor 
x_{i+1}, \ldots , 1 \tensor x_n$ where the terms $x_1 
\tensor 1, \ldots, x_i \tensor 1$ occur in the given 
natural order possibly interrupted by factors from the 
second set and similarly the terms $1 \tensor x_{i+1}, 
\ldots, 1 \tensor x_n$ occur in the given natural order 
possibly interrupted by factors from the first set. Such a 
product will occur in the expansion of 
 $$\rho(\sigma, (g_1, \ldots, g_n)) (x_{\sigma(1)} \tensor 
1 + 1 \tensor x_{\sigma(1)}) \ldots (x_{\sigma(n)} \tensor 
1 + 1 \tensor x_{\sigma(n)})$$ 
 whenever $\sigma \in S_n$ is a shuffle of $\{1, \ldots, 
i\}$ with $\{i+1, \ldots, n\}$, i.e. if $1 \leq j < k \leq 
i$ or $i+1 \leq j < k \leq n$ then $\sigma^{-1}(j) < 
\sigma^{-1}(k)$. 
 
 To evaluate such a product we have to interchange factors 
according to the rule given above, namely, 
 $$(1 \tensor x_m)(x_l \tensor 1) = \chi(g_m,g_l)(x_l 
\tensor x_m) = \chi(g_m,g_l) (x_l \tensor 1)(1 \tensor x_m) 
\hbox{ for } 1 \leq l \leq i < m \leq n.$$ 
 This rule has to be applied for every pair in the product 
(and in $\sigma$) in reverse position $j < k$ and $m := 
\sigma(j) > \sigma(k) = l$ thus producing a factor 
$\chi(g_{\sigma(j)},g_{\sigma(k)})$. So the total 
contribution of 
 $$\rho(\sigma, (g_1, \ldots, g_n)) (x_{\sigma(1)} \tensor 
1 + 1 \tensor x_{\sigma(1)}) \ldots (x_{\sigma(n)} \tensor 
1 + 1 \tensor x_{\sigma(n)})$$ 
 to $c_{1,i} \cdot x_1 \ldots x_i \tensor x_{i+1} \ldots 
x_n$ is 
 $$\begin{array}{r@{}l} 
 \rho(\sigma, (g_1, \ldots, g_n)&) \prod_{j < k, \sigma(j) 
> \sigma(k)} \chi(g_{\sigma(j)},g_{\sigma(k)}) \\ 
 &= \prod_{j < k, \sigma(j) > \sigma(k)} (\zeta^{- 
1}\chi(g_{\sigma(k)}, g_{\sigma(j)})) 
\chi(g_{\sigma(j)},g_{\sigma(k)})) = \zeta^t \\
 \end{array}
 $$ 
 where $t$ is the number of pairs in reverse position in 
$\sigma$. Here we have used that $(g_1, \dots, g_n)$ is a 
$\zeta$-family. 
 
 To determine the number $t$ of pairs in reverse position 
for a fixed $\sigma$ we observe that the ordering of $1, 
\ldots, i$ and of $i+1, \ldots, n$ must be preserved in the 
product. We count how many steps the factor $1 \tensor 
x_{i+1}$ has been moved to the left or how many terms from 
$x_1 \tensor 1, \ldots, x_i \tensor 1$ in the product are 
to the right of $1 \tensor x_{i+1}$ and call this number 
$\lambda_1$. Observe that $0 \leq \lambda_1 \leq i$. 
Similarly $\lambda_2$ denotes the number of terms from $x_1 
\tensor 1, \ldots, x_i \tensor 1$ in the product that are 
to the right of $1 \tensor x_{i+2}$. We have $0 \leq 
\lambda_2 \leq \lambda_1$. In a similar way we continue to 
define the numbers $\lambda_j$ with $0 \leq \lambda_{n-i} 
\leq \ldots \leq \lambda_2 \leq \lambda_1 \leq i$. The 
evaluation of the selected product then gives a term 
 $$\zeta^{\lambda_1 + \ldots + \lambda_{n-i}} x_1 \ldots 
x_i \tensor x_{i+1} \ldots x_n.$$ 
 
 To get the number of terms $\zeta^t x_1 \ldots x_i \tensor 
x_{i+1} \ldots x_n$ in $[x_1 \tensor 1 + 1 \tensor x_1, 
\ldots, x_n \tensor 1 + 1 \tensor x_n] $ we have to count 
all possibilities to represent $t = \lambda_1 + \ldots + 
\lambda_{n-i}$ with $0 \leq \lambda_{n-i} \leq \ldots \leq 
\lambda_2 \leq \lambda_1 \leq i$ or the number $p(i,n-i,t)$ 
of partitions of $t$ into at most $n-i$ parts each $\leq 
i$. Thus we can now determine the factor $c_{1,i}$ for 
$c_{1,i} \cdot x_1 \ldots x_i \tensor x_{i+1} \ldots x_n$ 
in the expansion of $[x_1 \tensor 1 + 1 \tensor x_1, 
\ldots, x_n \tensor 1 + 1 \tensor x_n]$ as 
 $$c_{1,i} = \sum_{t \geq 0} p(i,n-i,t) \zeta^t.$$ 
 By a theorem of Sylvester (\cite {A} Theorem 3.1) we have 
 $$\sum_{t \geq 0} p(i,n-i,t)q^t = 
 {{(1-q^n)(1-q^{n-1}) \ldots (1-q^{n-i+1})} 
 \over 
 {(1-q^i)(1-q^{i-1}) \ldots (1-q)}}$$ 
 hence $c_{1,i} = 0$ for $1 < i < n$ since $\zeta$ is a 
{\it primitive} $n$-th root of unity (see also \cite{T} p. 
2632). 
 
 So we have shown 
 $$\begin{array}{rl} 
 [x_1 \tensor 1 &+  1 \tensor x_1, \ldots, x_n \tensor 1 + 
1 \tensor x_n] = \\ 
 &= \sum_{\sigma \in S_n} \rho(\sigma, (g_1, \ldots, g_n)) 
x_{\sigma(1)} \ldots x_{\sigma(n)} \tensor 1 + 1 \tensor 
x_{\sigma(1)} \ldots x_{\sigma(n)} \\ 
 &= [x_1,\ldots, x_n] \tensor 1 + 1 \tensor [x_1, \ldots, 
x_n].\\
 \end{array}
 $$ 
 \end{pf}
 
 Now we will study Hopf algebras in the category $\M^{kG}$ 
and their primitive \footnote{This notion of primitivity 
is not related to the notion of primitivity for roots of 
unity.} elements. An element $x \in H$ is {\it primitive} 
if $\Delta(x) = x \tensor 1 + 1 \tensor x$. The set of 
primitive elements of an ordinary Hopf algebra forms a Lie 
algebra. This is not true for Hopf algebras in $\M^{kG}$. 
 
 Let $P(H)$ denote the set of primitive elements of $H$ and 
$P_g(H)$ the set of primitive elements of degree $g$ in 
$H$. Then 
 $$P(H) = \bigoplus_{g \in G} P_g(H)$$ 
 since for any $x = \sum_{g \in G}x_g$ we have 
 $$\sum_{g \in G} \Delta(x_g) = \Delta(x) = x \tensor 1 + 1 
\tensor x = \big( \sum_{g \in G} x_g \big) \tensor 1 + 1 
\tensor \big( \sum_{g \in G} x_g \big),$$ 
 so by comparing homogeneous terms we get $\Delta(x_g) = 
x_g \tensor 1 + 1\tensor x_g$. Thus the homogeneous 
components $x_g$ of $x$ are again primitive hence in 
$P_g(H)$ which shows $P(H) \subseteq \bigoplus_{g \in G} 
P_g(H)$. The converse inclusion is trivial. 
 
 \begin{thm} \label{prim_Lie_alg}
 Let $H$ be a $(G,\chi)$-Hopf algebra in $\M^{kG}$. Then 
for all primitive 
 $n$-th roots of unity $\zeta \not= 1$ and all 
$\zeta$-families $(g_1, \ldots, g_n)$ the following is a 
linear map 
 $$[ \ldots ]: P_{g_1}(H) \tensor \ldots \tensor P_{g_n}(H) 
\longrightarrow P_{g_1 + \ldots + g_n}(H)$$ 
 $$[x_1, \ldots, x_n] := \sum_{\sigma \in S_n} \rho(\sigma, 
(g_1, \ldots, g_n)) x_{\sigma(1)} \ldots x_{\sigma(n)}.$$ 
 \end{thm}
 
 \begin{pf} 
 It is clear that the degree of each product $x_{\sigma(1)} 
\ldots x_{\sigma(n)}$ is $g_1 + \ldots + g_n$. So we only 
have to show that $[x_1, \ldots, x_n]$ is primitive. But 
that is a consequence of Theorem \ref{main}. 
 \end{pf}
 
 Observe that we have special multiplications 
 $$[ \ldots ]: P_g(H) \tensor \ldots \tensor P_g(H) 
\longrightarrow P_{ng}(H)$$ for all $g \in G$ with $|g| 
\not= 1$ a primitive $n$-th root of unity. 
 
\begin{defn}
 Let $A$ be an algebra in $\M^{kG}$. A {\it derivation} 
from $A$ to $A$  of degree $g \in G$ is a family of linear 
maps $(d_h: A_h \longrightarrow A_{h+g} | h \in G )$ such 
that 
 $$d(ab) = d(a)b + \chi(g,h)ad(b)$$ for all $a \in A_h$, $b 
\in A_{h'}$, all $h,h' \in G$.
 \end{defn}
 
 It is clear that all derivations from $A$ to $A$ of all 
degrees form an object $\Der(A)$ in $\M^{kG}$ and that 
there is an operation $\Der(A) \tensor A \longrightarrow 
A$. 
 
 \begin{cor}
 $\Der(A)$ is a $(G,\chi)$-Lie algebra.
 \end{cor}
 
 \begin{pf} 
 Let $m$ denote the multiplication of $A$. An endomorphism 
$x: A \longrightarrow A$ of degree $g \in G$, i.e. $x_h: 
A_h \longrightarrow A_{h+g}$ for all $h \in G$ is a 
derivation iff $m(x \tensor 1 + 1 \tensor x) = xm$ where 
$(x \tensor y)(a \tensor b) = \chi(\deg(y), \deg(a)) x(a) 
\tensor y(b)$ for homogeneous elements $a$ and $b$ in $A$. 
 
 To show that $\Der(A)$ is a Lie algebra it suffices to 
show that it is closed under Lie multiplication since it is 
a subobject of $\End(A)$, the inner endomorphism object of 
$A$ which is known to be an algebra in the category 
$\M^{kG}$. Let $\zeta$ be a primitive $n$-th root of unity 
and let $(g_1, \ldots, g_n)$ be a $\zeta$-family. Then 
 $$\begin{array}{rl} 
 & m([x_1, \ldots, x_n] \tensor 1 + 1 \tensor [x_1, \dots, 
x_n]) \\ 
 &= m[x_1 \tensor 1 + 1 \tensor x_1, \ldots , x_n \tensor 1 
+ 1 \tensor x_n] \\ 
 &= m \sum_{\sigma \in S_n} \rho(\sigma,(g_1, \ldots, g_n)) 
(x_{\sigma(1)} \tensor 1 + 1 \tensor x_{\sigma(1)})\cdot 
\ldots \cdot (x_{\sigma(n)} \tensor 1 + 1 \tensor 
x_{\sigma(n)}) \\ 
 &= \sum_{\sigma \in S_n} \rho(\sigma,(g_1, \ldots , g_n)) 
m(x_{\sigma(1)} \tensor 1 + 1 \tensor x_{\sigma(1)}) \cdot 
\ldots \cdot (x_{\sigma(n)} \tensor 1 + 1 \tensor 
x_{\sigma(n)}) \\ 
 &= \sum_{\sigma \in S_n} \rho(\sigma,(g_1, \ldots, g_n)) 
x_{\sigma(1)} \cdot \ldots \cdot x_{\sigma(n)} \cdot m \\ 
 &= [x_1, \ldots, x_n] m \\
 \end{array}
 $$ 
 for all derivations $x_1, \ldots, x_n$ of degrees $g_1, 
\ldots, g_n$ respectively. Hence $[x_1, \ldots, x_n]$ again 
is a derivation. 
 \end{pf}

 \section{Lie Algebras and Universal Enveloping 
Algebras} 
 
 \begin{defn} \label{Lie_alg}
 An object $P = \bigoplus P_g \in \M^{kG}$ together with 
operations 
 $$[ \ldots ]: P_{g_1} \tensor \ldots \tensor P_{g_n} 
\longrightarrow P_{g_1 + \ldots + g_n}$$ 
 for all $n \in \N$, all primitive $n$-th roots of unity 
$\zeta$, and all $\zeta$-families $(g_1, \ldots, g_n)$, is 
called a $(G,\chi)${\it -Lie algebra} if the following 
identities hold: 
 \begin{enumerate}
 \item for all primitive $n$-th roots of unity $\zeta$, 
all $\zeta$-families $(g_1, \ldots, g_n)$, all $\sigma \in 
S_n$, and all $x_i \in P_{g_i}$ 
 $$[x_1, \ldots, x_n] = \rho(\sigma,(g_1, \ldots, g_n)) 
[x_{\sigma(1)}, \ldots, x_{\sigma(n)}],$$ 
 \item 
 for all primitive $n$-th roots of unity $\zeta$, all 
$\zeta$-families $(g_1, \ldots, g_{n+1})$, and for all $x_i 
\in P_{g_i}$ 
 $$\sum_{i=1}^{n+1} \rho((i \ldots 1), (g_1, \ldots, 
g_{n+1})) [x_i,[x_1,\ldots,\hat x_i, \ldots, x_{n+1}]] = 
0,$$ 
 \item 
 for all primitive $n$-th roots of unity $\zeta$, all 
$\zeta$-families $(g_1, \ldots, g_n)$, all $h \in G$ such 
that all $(h,g_i)$ are $(-1)$-families, all $y_i \in 
P_{g_i}$, and all $x \in P_h$ 
 $$[x,[y_1, y_2, \ldots, y_n]] = \sum_{i=1}^n \big( 
\prod_{j=1}^{i-1} \chi(g_j,h) \big) [y_1, \ldots, [x,y_i], 
\ldots, y_n].$$ 
 \end{enumerate}
 \end{defn}
 
 The $(G,\chi)$-Lie algebras form a category in a 
straightforward way and the construction given in section 3 
defines a functor $P : (G,\chi)\Hopf \longrightarrow 
(G,\chi)\Lie$ as can be easily verified. 
 
 The Theorems \ref{prim_Lie_alg}, \ref{symm}, and 
\ref{jacobi} show that this notion of generalized Lie 
algebras passes its test on the set of primitive elements 
of a Hopf algebra: 
 
 \begin{cor}
 Let $H$ be a $(G,\chi)$-Hopf algebra in $\M^{kG}$. Then 
the set of primitive elements $P(H)$ forms a $(G,\chi)$-Lie 
algebra. 
 \end{cor}
 
 \begin{lma}

 Let $A$ be a $(G,\chi)$-algebra in $\M^{kG}$. Then $A$ 
carries the structure of a $(G,\chi)$-Lie algebra $A^L$ 
with the operations 
 $$[x_1, \ldots, x_n] := \sum_{\sigma \in S_n} \rho(\sigma, 
(g_1, \ldots, g_n)) x_{\sigma(1)} \cdot \ldots \cdot 
x_{\sigma(n)}$$ 
 for all (roots of unity) $\zeta \in k^*$, all 
$\zeta$-families $(g_1, \ldots, g_n)$, and all $x_i \in 
A_{g_i}$. 
 \end{lma}
 
 \begin{pf} 
 This is a rephrasing of Theorems \ref{symm} and \ref{jacobi}. 
 \end{pf}
 
 This lemma defines a functor $(G,\chi)\Alg \ni A \mapsto 
A^L \in (G,\chi)\Lie$. 
 
 \begin{thm}
 Let $P$ be a $(G,\chi)$-Lie algebra. Then there is a 
universal associative enveloping algebra $U(P)$ in 
$\M^{kG}$. 
 \end{thm}
 
 \begin{pf} 
 We define $U(P) := T(P)/I$ where $T(P)$ is the tensor 
algebra and $I := ([x_1, \ldots, x_n] - \sum_{\sigma \in 
S_n} \rho(\sigma, (g_1, \ldots, g_n)) x_{\sigma(1)} \cdot 
\ldots \cdot x_{\sigma(n)}) $ is the ideal generated by all 
terms formed for all (roots of unity) $\zeta \in k^*$, all 
$\zeta$-families $(g_1, \ldots, g_n)$, and all $x_i \in 
P_{g_i}$. 
 
 The tensor algebra is constructed in $\M^{kG}$ with the 
natural grading. Since the relations are $G$-homogeneous, 
the algebra $U(P)$ is also in $\M^{kG}$. 
 
 By the universal property of the tensor algebra any 
$(G,\chi)$-Lie homomorphism $f: P \longrightarrow A^L$ 
extends uniquely to a $(G,\chi)$-algebra homomorphism $g: 
U(P) \longrightarrow A$, so that $U(P)$ is the universal 
associative algebra generated by $P$. 
 \end{pf}
 
 \begin{cor}
 The functor $\z^L: (G,\chi)\Alg \longrightarrow 
(G,\chi)\Lie$ has the left adjoint functor $U: (G,\chi)\Lie 
\longrightarrow (G,\chi)\Alg$. 
 \end{cor}
 
 Observe that the identities for $(G,\chi)$-Lie algebras 
play no special role in the construction of the universal 
enveloping algebra. Its construction depends only on the 
given operations, which must be homogeneous so that the 
residue class algebra is $G$-graded again. 
 
 If we consider the set $P(H)$ of primitive elements of an 
arbitrary $(G,\chi)$-Hopf algebra $H$ then there are many 
homogeneous partially defined common operations in all 
$P(H)$. We used the most basic ones in our generalization 
of Lie algebras. We would like to know if the 
bracket-operations and identities given in the definition 
of a $(G,\chi)$-Lie algebra are a generating set for all 
possible operations and relations for the $G$-graded module 
$P(H)$ for every Hopf algebra in $\M^{kG}$, every abelian 
group $G$, every bicharacter $\chi$ and every integral 
domain $k$. This should be true for the identities if all 
maps from $(G,\chi)$-Lie algebras into their universal 
enveloping algebras are injective. To prove this, some kind 
of generalization of the Poincar\'e-Birkhoff-Witt theorem 
is needed. 
 
 Now we show how to construct Hopf algebras in the braided 
category of $G$-graded vector spaces using universal 
enveloping algebras. 

 \begin{thm}
 Let $P$ be a $(G,\chi)$-Lie algebra. Then the universal 
enveloping algebra $U(P)$ is a $(G,\chi)$-Hopf algebra. 
 \end{thm}
 
 \begin{pf} 
 We define a homomorphism $\delta: P \longrightarrow (U(P) 
\tensor U(P))^L$ in $\M^{kG}$ by $\delta(x) := \x \tensor 1 
+ 1 \tensor \x$ where $\x$ is the canonical image of $x \in 
P$ in $U(P)$. Then we have 
 $$\begin{array}{rl} 
 [\delta(x_1),& \ldots, \delta(x_n)] \\ 
 &= [\x_1 \tensor 1 + 1 \tensor \x_1, \ldots, \x_n \tensor 
1 + 1 \tensor \x_n] \\ 
 &= [\x_1, \ldots, \x_n] \tensor 1 + 1 \tensor [\x_1, 
\ldots, \x_n] \quad \hbox{(by Theorem \ref{main})} \\ 
 &= \overline{[x_1, \ldots,x_n]} \tensor 1 + 1 \tensor 
\overline{[x_1, \ldots,x_n]} \quad \hbox{(by definition of 
} U(P)\ ) \\ 
 &= \delta([x_1, \ldots,x_n]). \\
 \end{array}
 $$ 
 Hence $\delta: P \longrightarrow (U(P) \tensor U(P))^L$ is 
a $(G,\chi)$-Lie homomorphism that factors uniquely through 
$\Delta: U(P) \longrightarrow U(P) \tensor U(P)$, an 
algebra homomorphism. So for all $x \in P$ we have 
$\Delta(\x) = \x \tensor 1 + 1 \tensor \x$. 
 
 Since $(\Delta \tensor 1)\delta(x) = \x \tensor 1 \tensor 
1 + 1 \tensor \x \tensor 1 + 1 \tensor 1 \tensor \x = (1 
\tensor \Delta)\delta(x)$ and since $(\Delta \tensor 
1)\delta = (1 \tensor \Delta)\delta$ is a $(G,\chi)$-Lie 
homomorphism, it factors through a unique algebra 
homomorphism $(\Delta \tensor 1)\Delta = (1 \tensor 
\Delta)\Delta$ so that $\Delta: U(P) \longrightarrow U(P) 
\tensor U(P)$ is coassociative. 
 
 The counit $\varepsilon: U(P) \longrightarrow k$ is 
defined by the zero morphism $0: P \longrightarrow k$. Thus 
$U(P)$ is a bialgebra. 
 
 The definition of the antipode is somewhat more 
complicated. We consider $U(P)^{op}$ the opposite algebra 
of $U(P)$ where the multiplication in $U(P)^{op}$ is 
defined by $x_2 \circ x_1 = \chi(g_2,g_1)x_1x_2$ for 
homogeneous elements $x_i \in U(P)_{g_i}$ hence $x_1x_2 = 
\chi(g_2,g_1)^{-1}x_2 \circ x_1$. 
 
 Define a homomorphism $S: P \longrightarrow U(P)^{op}$ by 
$S(x) := -\x$. To show that this is a $(G,\chi)$-Lie 
homomorphism let $\zeta$ be a primitive $n$-th root of 
unity, $(g_1, \ldots, g_n)$ be a $\zeta$-family, and $x_i 
\in P_{g_i}$. We first consider the permutation $\tau \in 
S_n$ with $\tau(i) = n+1-i$. Then $\tau^2 = 1$ and 
 $$\begin{array}{rl} 
 \rho(\tau, (g_1, \ldots, g_n)) 
 &= \prod_{i < j} (\zeta^{-1} \chi(g_{n+1-j}, g_{n+1-i})) = 
\prod_{i < j} (\zeta \chi(g_{j}, g_{i})^{-1}) \\ 
 &= \zeta^{n(n-1)/2} \prod_{i < j} \chi(g_j,g_i)^{-1} = (- 
1)^{n-1} \prod_{i < j} \chi(g_j,g_i)^{-1}. \\
 \end{array}
 $$ 
 Using Lemma \ref{rho_is_homo} we get for all $\sigma \in S_n$ 
 $$\begin{array}{r@{}l} 
 \rho(\sigma, &(g_1, \ldots, g_n)) 
 \prod_{i < j} \chi(g_{\sigma(j)},g_{\sigma(i)})^{-1} \\ 
 &= (-1)^{n-1} \rho(\sigma,(g_1, \ldots, g_n)) 
 \rho(\tau,(g_{\sigma(1)}, \ldots, g_{\sigma(n)})) \\ 
 &= (-1)^{n-1} \rho(\sigma\tau,(g_1, \ldots, g_n)). \\
 \end{array}
 $$ 
 Thus we have 
 $$\begin{array}{r@{}l} 
 S([x_1,\ldots,x_n & ]) = -\overline{[x_1, \ldots, x_n]} = - 
[\x_1,\ldots,\x_n] \\ 
 &= -\sum_{\sigma \in S_n} \rho(\sigma,(g_1, \ldots, g_n)) 
    \x_{\sigma(1)} \ldots \x_{\sigma(n)} \\ 
 &= -\sum_{\sigma \in S_n}\rho(\sigma,(g_1, \ldots, g_n)) 
    \big( \prod_{i < j} 
\chi(g_{\sigma(j)},g_{\sigma(i)})^{-1} \big) 
    \x_{\sigma(n)} \circ \ldots \circ \x_{\sigma(1)} \\ 
 &= (-1)^n\sum_{\sigma \in S_n}\rho(\sigma\tau,(g_1, 
\ldots, g_n)) 
    \x_{\sigma(n)} \circ \ldots \circ \x_{\sigma(1)} \\ 
 &= (-1)^n\sum_{\sigma \in S_n}\rho(\sigma,(g_1, \ldots, 
g_n)) 
    \x_{\sigma(1)} \circ \ldots \circ \x_{\sigma(n)} \\ 
 &= [-\x_1, \ldots, -\x_n]^{op} = [S(x_1) , \ldots, 
S(x_n)]^{op}. \\
 \end{array}
 $$ 
 Thus $S$ can be extended to an algebra homomorphism $S: 
U(P) \longrightarrow U(P)^{op}$. 
 
 We show now that $S$ is the antipode for the bialgebra 
$U(P)$. Let $x \in P_h$, $x_i \in P_{g_i}$. We observe 
first that 
 $$S(\x\x_1\ldots \x_n) = S(\x) \circ S(\x_1\ldots \x_n) 
 = \chi(h,g_1 + \ldots + g_n) S(\x_1 \ldots \x_n)S(\x)$$ 
 since $S$ is a $(G,\chi)$-algebra homomorphism. We now 
prove $\nabla(1 \tensor S)\Delta(a) = \varepsilon(a)$ for 
homogeneous terms $a = 1$, $a = \x$ resp. $a = \x\x_1 
\ldots \x_n$ in $U(P)$ by induction. 
 $$\nabla(1 \tensor S)\Delta(1) = 1S(1) = 1 = \varepsilon 
(1), $$ 
 $$\nabla(1 \tensor S)\Delta(\x) = \x + S(\x) = 0 = 
\varepsilon(\x).$$ 
 Since we can write $\Delta(a) = \sum_i a_{i,1} \tensor 
a_{i,2} $ with each of the terms $a_{i,j}$ homogeneous and 
$\deg(a_{i,1}) + \deg(a_{i,2}) = \deg(a)$ ($U(P)$ is a 
$kG$-comodule coalgebra) we get 
 $$\begin{array}{r@{\ }l} 
 \nabla(1 \tensor S)\Delta(xa) 
 &= \nabla(1 \tensor S) ( x \tensor 1 + 1 \tensor x)(\sum 
a_{i,1} \tensor a_{i,2}) \\ 
 &= \nabla(1 \tensor S) \big( \sum xa_{i,1} \tensor a_{i,2} 
    + \sum \chi(h,\deg(a_{i,1})) a_{i,1} \tensor xa_{i,2} 
\big) \\ 
 &= \sum xa_{i,1} S(a_{i,2}) 
    + \sum \chi(h,\deg(a_{i,1}))a_{i,1} S(xa_{i,2})  \\ 
 &= 0 + \sum \chi(h,\deg(a_{i,1})) \chi(h,\deg(a_{i,2})) 
    a_{i,1} S(a_{i,2}) S(x)  \\ 
 &= \chi(h,\deg(a)) \sum a_{i,1} S(a_{i,2}) S(x) \\ 
 &= 0 = \varepsilon(xa). \\
 \end{array}
 $$ 
 Analogously we get $\nabla(S \tensor 1)\Delta = 
\varepsilon$, so that $S$ is an antipode and $U(P)$ is a 
$(G,\chi)$-Hopf algebra. 
 \end{pf}
 
 We come to an interesting consequence for a universal 
enveloping algebra. There is a new multiplication on $U(P)$ 
defined by $a \cdot b := \chi(h,g)\chi(g,h) ab$ for $a,b$ 
homogeneous with $\deg(a) = h$, $\deg(b) = g$. This algebra 
is $U(P)^{op\ op}$ which in general will be different from 
$U(P)$. The map $P \ni x \mapsto \bar x \in U(P)^{op\ op}$ 
is a Lie homomorphism (similar to the proof that $S$ was a 
Lie homomorphism). Thus it induces a homomorphism $U(P) 
\longrightarrow U(P)^{op\ op}$. Since $U(P)^{op\ op}$ is 
also generated by $P$ it is easy to see, that it is also a 
universal enveloping algebra hence $U(P)$ and $U(P)^{op\ 
op}$ are isomorphic under the given homomorphism. 
 
 \begin{cor}
 The functor $U: (G,\chi)\Lie \longrightarrow 
(G,\chi)\Hopf$ is left adjoint to the functor $P: 
(G,\chi)\Hopf \longrightarrow (G,\chi)\Lie$. 
 \end{cor}
 
 Now we give some examples of braided monoidal categories 
of $G$-graded modules together with the associated 
$(G,\chi)$-Lie structures induced by associative algebras. 
Let $k$ be an algebraically closed field of characteristic 
zero. 
 
  \begin{eg} \label{brackets}
 (1) Let $G = \{0\}$. Then we have $\chi(0,0) = 1$ and 
hence $\chi(0,0)\chi(0,0) = 1 =(-1)^2$ so that 
$(0,0)$ is a $(-1)$-family. Thus we have a bracket 
multiplication 
 $$[.,.]: P_0 \tensor P_0 \longrightarrow P_0,$$ 
 which for associative algebras $A$ is defined by 
 $$[x,y] = xy - yx.$$ 
 This is the example of ordinary Lie algebras. The Lie 
identities of Definition \ref{Lie_alg} reduce to 
 $$[x_1,x_2] = \rho((2,1),(0,0))[x_2,x_1] = -[x_2,x_1],$$ 
 $$\begin{array}{rl} 
 &\sum_{i=1}^3 \rho((i \ldots 1), (0,0,0)) 
[x_i,[x_1,\ldots,\hat x_i, \ldots, x_3]] = \\ 
 &= [x_1,[x_2,x_3]] + \rho((2, 1), (0,0,0)) [x_2,[x_1,x_3]] 
+ \rho((3,2,1), (0,0,0)) [x_3,[x_1,x_2]] \\ 
 &= [x_1,[x_2,x_3]] - [x_2,[x_1,x_3]] + [x_3,[x_1,x_2]] = 
0, \\
 \end{array}
 $$ and 
 $$[x,[y_1, y_2]] = [[x,y_1], y_2] + \chi(0,0) [y_1, 
[x,y_2]] = [[x,y_1], y_2] + [y_1, [x,y_2]].$$ 
 
 (2) Let $G = \Z/2\Z = \{0,1\}$ with the bicharacter 
$\chi(i,j) = (-1)^{ij}$. We have $\chi(0,i) = 1$ and hence 
$\chi(0,1)\chi(1,0) = 1 = (-1)^2$ so the following are 
$(-1)$-families: $(0,0)$, $(0,1)$, and $(1,0)$ that induce 
bracket multiplications 
 $$[.,.]: P_0 \tensor P_0 \longrightarrow P_0,$$ 
 $$[.,.]: P_0 \tensor P_1 \longrightarrow P_1,$$ 
 $$[.,.]: P_1 \tensor P_0 \longrightarrow P_1,$$ 
 which for associative algebras $A$ are defined by 
 $$[x,y] = xy - yx.$$ 
 Furthermore we have $\chi(1,1) = -1$ and hence 
$\chi(1,1)\chi(1,1) = 1 = (-1)^2$ so $(1,1)$ is another $(- 
1)$-family that induces a bracket multiplication 
 $$[.,.]: P_1 \tensor P_1 \longrightarrow P_0,$$ with 
 $$[x,y] = xy + \rho((2,1),(1,1))yx = xy + yx.$$ 
 This is the example of Lie super algebras. The Lie 
identities of Definition \ref{Lie_alg} reduce to 
 $$[x_1,x_2] = \rho((2,1),(\deg(x_1),\deg(x_2)))[x_2,x_1] = 
-(-1)^{\deg(x_1)\deg(x_2)}[x_2,x_1],$$ 
 $$\begin{array}{r@{\ }l} 
 [x_1,[x_2,x_3]] - &(-1)^{\deg(x_1)\deg(x_2)} 
[x_2,[x_1,x_3]] \\ 
 &\ \ + (-1)^{(\deg(x_1) + 
\deg(x_2))\deg(x_3)}[x_3,[x_1,x_2]] = 0, \\
 \end{array}
 $$ 
 $$[x,[y_1,y_2]] = [[x,y_1],y_2] + (-1)^{\deg(x)\deg(y_1)} 
[y_1,[x,y_2]].$$ 
 
 (3) Let $G$ be an arbitrary abelian group with a 
bicharacter $\chi$ such that $\chi(g_1,g_2) = 
\chi(g_2,g_1)^{-1}$ for all $g_1,g_2 \in G$. Then we have 
$\chi(g_1,g_2)\chi(g_2,g_1) = 1 = (-1)^2$. This defines $(- 
1)$-families $(g_1,g_2)$ together with bracket operations 
 $$[.,.]: P_{g_1} \tensor P_{g_2} \longrightarrow 
P_{g_1+g_2},$$ 
 $$[x,y] = xy + \rho((2,1),(g_1,g_2))yx = xy - 
\chi(g_1,g_2)yx.$$ 
 This is the example of Lie color algebras for an  abelian 
group. The Lie identities of Definition \ref{Lie_alg} reduce to 
 $$[x_1,x_2] = \rho((2,1),(\deg(x_1),\deg(x_2)))[x_2,x_1] = 
-\chi(\deg(x_1),\deg(x_2)) [x_2,x_1],$$ 
 $$\begin{array}{rl} 
 &\ \ [x_1,[x_2,x_3]] - \chi(\deg(x_1),\deg(x_2)) 
[x_2,[x_1,x_3]] \\ 
 &+ \chi(\deg(x_1) + \deg(x_2),\deg(x_3)) [x_3,[x_1,x_2]] = 
0, \\
 \end{array}
 $$ 
 $$[x,[y_1,y_2]] = [[x,y_1],y_2] + \chi(\deg(x),\deg(y_1)) 
[y_1,[x,y_2]].$$ 
 
 (4) Let $G = \Z/3\Z = \{0,1,2\}$ with the bicharacter 
$\chi(i,j) = \zeta^{ij}$ where $\zeta$ is a primitive 3rd 
root of unity. Then $\chi(0,i)\chi(i,0) = 1 = (-1)^2$ so 
that we get $(-1)$-families $(0,i)$ with 
 $$[.,.]: P_0 \tensor P_i \longrightarrow P_i,$$ 
 $$[.,.]: P_i \tensor P_0 \longrightarrow P_i,$$ 
 $$[x,y] = xy - yx.$$ 
 Furthermore we have $\chi(i,j)\chi(j,i) = (\zeta^{ij})^2$ 
with a primitive 3rd root of unity $\zeta^{ij}$ for all 
$i,j \not= 0$. This gives $\zeta$-families $(1,1,1)$ and 
$(2,2,2)$ and no $\zeta^2$-family. The associated Lie 
structure is 
 $$[ \ldots ]: P_i \tensor P_i \tensor P_i \longrightarrow 
P_0$$ with 
 $$[x_1,x_2,x_3] = \sum_{\sigma \in S_3} 
x_{\sigma(1)}x_{\sigma(2)}x_{\sigma(3)}$$ 
 since $\rho(\sigma,(i,i,i)) = \prod (\zeta^{-1}\chi(i,i)) 
= 1$ for all $\sigma \in S_3$.  The Lie identities of 
Definition \ref{Lie_alg} for these ternary brackets reduce to 
 $$[x_1,x_2,x_3] = [x_{\sigma(1)}, x_{\sigma(2)}, 
x_{\sigma(3)}],$$ 
 $$[x_1,[x_2,x_3,x_4]] + [x_2,[x_1,x_3,x_4]] + 
[x_3,[x_1,x_2,x_4]] + [x_4,[x_1,x_2,x_3]] = 0,$$ 
 for $x_1, x_2, x_3, x_4 \in P_i$, $i=1$ or $i=2$. Here we 
also find a first example where the two Jacobi identities 
mean different things. The second Jacobi identity for $x 
\in P_0$ and $y_1, y_2, y_3 \in P_i$, $i = 1$ or $i = 2$ 
reduces to 
 $$[x,[y_1,y_2,y_3]] = [[x,y_1],y_2,y_3] + 
[y_1,[x,y_2],y_3] + [y_1,y_2,[x,y_3]].$$ 
 Furthermore there are $(-\zeta)$-families $(1,1,1,1,1,1)$ 
and $(2,2,2,2,2,2)$ for the primitive 6-th root of unity $-
\zeta$ which give maps $\tensor^6 P_1 \to P_0$ and 
$\tensor^6 P_2 \to P_0$. The associated Lie multiplication 
of an associative algebra is
 $$[x_1,\ldots,x_6] = \sum_{\sigma \in S_6} 
\sgn(\sigma) x_{\sigma(1)}\ldots x_{\sigma(6)}$$ 
 since $\rho(\sigma,(i,i,i,i,i,i)) = \prod (-\zeta^{-1}\chi(i,i)) 
= \sgn(\sigma)$ for all $\sigma \in S_6$.
 
 (5) Let $G = C_n \times \ldots \times C_n$ ($r$-times) 
with generators $t_1, \ldots, t_r$. Since later we will 
take the biproduct of a $(G,\chi)$-Hopf algebra with $kG$, 
we will write $G$ multiplicatively in this example. Let 
$\chi$ be given by $\chi(t_i,t_j) = \zeta$, a primitive 
$n$-th root of unity. Then $(t_1, \ldots, t_1)$ is a 
$\zeta$-family. A $(G,\chi)$-Lie algebra will have a 
bracket operation 
 $$[.,.]: \tensor^n P_{t_1}  \longrightarrow P_0$$ 
 with 
 $$[x_1,\ldots,x_n] = \sum_{\sigma \in S_n} x_{\sigma(1)} 
\ldots x_{\sigma(n)}.$$ 
 
 (6) The last example of this kind comes with $G = \Z/3\Z 
\times \Z/3\Z$ and generators $g_1, g_2$ each of order 3. 
Define $\chi$ by $\chi(g_1,g_1) := \zeta$, a primitive 3rd 
root of unity, $\chi(g_1,g_2) := \zeta^2$, $\chi(g_2,g_1) = 
1$, and $\chi(g_2,g_2) = \zeta$. Then there are several 
$\zeta$-families, among others $(g_1, g_1, g_2)$ and $(g_1, 
g_2, g_2)$. They define brackets 
 $$[ \ldots ]: P_{g_1} \tensor P_{g_1} \tensor P_{g_2} 
\longrightarrow P_{2g_1 + g_2},$$ 
 $$[ \ldots ]: P_{g_1} \tensor P_{g_2} \tensor P_{g_2} 
\longrightarrow P_{g_1 + 2g_2},$$ with 
 $$\begin{array}{rl} 
 [x_1,x_2,x_3] &= \sum_{\sigma \in S_3} 
\rho(\sigma,(g_1,g_1,g_2)) 
x_{\sigma(1)}x_{\sigma(2)}x_{\sigma(3)} \\ 
 &= x_1x_2x_3 + x_2x_1x_3 + \zeta x_1x_3x_2 + \zeta 
x_2x_3x_1 + \zeta^2 x_3x_1x_2 + \zeta^2 x_3x_2x_1 \\
 \end{array}
 $$ 
 resp. 
 $$\begin{array}{rl} 
 [x_1,x_2,x_3] &= \sum_{\sigma \in S_3} 
\rho(\sigma,(g_1,g_1,g_2)) 
x_{\sigma(1)}x_{\sigma(2)}x_{\sigma(3)} \\ 
 &= x_1x_2x_3 + x_1x_3x_2 + \zeta x_2x_1x_3 + \zeta 
x_3x_1x_2 + \zeta^2 x_2x_3x_1 + \zeta^2 x_3x_2x_1.\\
 \end{array}$$ 
 \end{eg}
 
 We close with some examples of $(G,\chi)$-Lie algebras and 
of $(G,\chi)$-Hopf algebras generated by them. We also 
examine some of the biproduct Hopf algebras one always 
obtains from $(G,\chi)$-Hopf algebras. If $H$ is a 
$(G,\chi)$-Hopf algebra then $H \tensor kG$ is a Hopf 
algebra by 
 $$(x \tensor g)(y \tensor h) = \chi(g, \deg(y)) xy \tensor 
gh,$$ 
 $$\Delta(1 \tensor g) = (1 \tensor g) \tensor (1 \tensor 
g),$$ 
 and 
 $$\Delta(x \tensor 1) = \sum_{a + b = c} (y_{a,i} \tensor 
b) \tensor (z_{b,i} \tensor 1),$$ 
 where $\Delta(x) = \sum_{a+b = c}\sum_i y_{a,i} \tensor 
z_{b,i}$ with $\deg(x) = c$, $\deg(y_{a,i}) = a$, and 
$\deg(z_{b,i}) = b$ (see \cite {FM} Corollary 3.5). 
 
 Most of the $(G,\chi)$-Hopf algebras are 
infinite-dimensional. In fact, the only finite-dimensional 
Hopf algebras generated by their primitives, we know, are 
given in the following example under (1) and (2). 
 
 \begin{eg}
 (1) $P = P_1 = kx$ with $[x,x] = 0$ defines a 
``commutative'' $(C_2,\chi)$-Lie algebra with $\chi$ as in 
Example \ref{brackets} (2). It generates the $(C_2,\chi)$-Hopf 
algebra $H = k[x]/(x^2)$, the universal enveloping algebra 
of $P$. The  biproduct $H \star\ kC_2$ is the well known 
smallest noncommutative noncocommutative Hopf algebra. 
(\cite {FM} Example 3.9) 
 
 (2) Let $(G, \chi)$ be as in Example \ref{brackets} (5). $P = P_{t_1} 
= kx$ with $[x, \ldots, x] = 0$ defines a ``commutative'' 
$(G,\chi)$-Lie algebra with $\chi(t_1,t_1) = \zeta$, a 
primitive $n$-th root of unity. It generates the 
$(G,\chi)$-Hopf algebra $H = k[x]/(x^n)$, given in the 
introduction (for $r = 1$). The biproduct $U(P) \star\ kG$ 
is the free algebra generated by the elements $x$, $t_1, 
\ldots, t_r$ subject to the following relations $t_i^n = 
1$, $t_it_j = t_jt_i$, and $x^n = 0$. There is one more set 
of relations that is obtained from the multiplicative rule 
for the biproduct $(1 \tensor t_i) (x \tensor 1) = 
\chi(t_i,t_1)x \tensor t_i$ or -- identifying $x \tensor 1$ 
with $x$ and $1 \tensor t_i$ with $t_i$ -- the relations 
$t_ix = \zeta xt_i$. Hence the biproducts form the family 
of noncommutative noncocommutative Hopf algebras given in 
\cite{T}. 
 
 (3) $P = \bigoplus_{i=0}^2 P_i$, $P_0 = kz$, $P_1 = kx 
\oplus ky$, $P_2 = 0$ with $[x,x,x] = [y,y,y] = [x,y,y] = 
0$,  $[x,x,y] = z$, and $[z,x] = [z,y] = 0$  define a 
$(C_3,\chi)$-Lie algebra with $\chi(\bar 1,\bar 1) = 
\zeta$, a primitive 3rd root of unity. It generates the 
$(C_3,\chi)$-Hopf algebra $H = k\langle 
x,y\rangle/(x^3,y^3,xy^2 + yxy + y^2x)$ (we have $z = 
2(x^2y + xyx + yx^2)$), the universal enveloping algebra of 
$P$. The biproduct $U(P) \star kC_3$ has generators $x, y, 
t$ with the relations $x^3 = y^3 = 0$, $t^3 = 1$, $xy^2 + 
yxy + y^2x = 0$, $xt = \zeta tx$, and $yt = \zeta ty$. 
 
 (4) Let $G = C_3$, and $\chi$ be as before. Then $P = 
\bigoplus_{i=0}^2 P_i$, $P_0 = ky$, $P_1 = kx$, $P_2 = 0$ 
with $[x,x,x,x,x,x] = 0$, $[x,x,x] = y$, and $[y,x] = 0$ 
defines a $(C_3,\chi)$-Lie algebra. It generates the 
$(C_3,\chi)$-Hopf algebra $H = k[x,y]/(y - 6 x^3) = k[x]$, 
the universal enveloping algebra of $P$. 
 
 (5)For $G$, $\chi$ as before let $P_0 = ky$, $P_1 = kx$, 
and $P_2 = 0$ with $[x,x,x,x,x,x] = 0$, $[x,x,x] = 0$ and 
$[y,x] = x$. Then $P$ is a $(C_3,\chi)$-Lie algebra. It 
generates the $(C_3,\chi)$-Hopf algebra $H = k\langle x,y 
\rangle/(x^3, x - yx + xy)$, the universal enveloping 
algebra of $P$. 
 \end{eg}
 

\end{document}